\title{On the cover time of planar graphs}
\author{Johan Jonasson\thanks{Chalmers University of Technology}
\and Oded Schramm\thanks{Microsoft Research}} 
\date{May 17, 2000}

\documentclass[12pt,fleqn,twoside]{article}
\usepackage{latexsym}
\usepackage{amsmath}
\usepackage{amsfonts}
\numberwithin{equation}{section}
\newtheorem{prop}{\sc Proposition}[section]     
\newtheorem{lemma}[prop]{\sc Lemma}
\newtheorem{theorem}[prop]{\sc Theorem}

\newcommand{\E}{{\mathbf E}}

\newcommand{\Z}{{\mathbb Z}}

\newcommand{\C}{{\mathbb C}}

\newcommand{\RR}{{\mathbb R}}
\newcommand{\N}{{\mathbb N}}

\def\dist{\hbox{\rm dist}}
\def\diam{\hbox{\rm diam}}
\def\area{\hbox{\rm area}}
\def\Re{\hbox{\rm Re}}

\begin{document}

\maketitle

\begin{abstract}
The cover time of a finite connected graph is the expected number
of steps needed for a simple random walk on the graph to visit
all the vertices.
It is known that the cover time on any $n$-vertex, connected graph is at
least
$\bigl(1+o(1)\bigr)n\log n$ and at most $\bigl(1+o(1)\bigr)\frac{4}{27}n^3$.
This paper proves that
for bounded-degree planar graphs the cover time
is at least $c n(\log n)^2$, and at most $6n^2$,
where $c$ is a positive constant depending only on the maximal degree
of the graph.
The lower bound is established via use of circle packings.
\end{abstract}

\section{Introduction}
Let $G=(V,E)$ be a finite, connected, $n$-vertex
graph and let $\{X_{k}\}_{k=0}^{\infty}$
be a simple random walk on $G$.
For each $v \in V$, set $T_{v} = \min\{k\in\N: X_k = v\}$ and let 
$C = \max_{v\in V}T_{v}$ be the {\em cover time}.
We are primarily interested in the expected cover
time $\E_{v} C$, where $\E_{v}$ denotes
expectation
with respect to the probability measure of the random walk starting at
$X_0 = v$.
In words, $\E_{v}C$ is the expected time taken for the random walk
starting at $v$ to visit every
vertex of the graph.

Over the last decade or so, much work has been devoted to finding
the expected cover time for different graphs and to giving general 
upper and lower bounds of the cover time.
For an introduction, we refer the reader
to the draft book by Aldous and Fill [\ref{AF}],
in particular to Chapters 3, 5 and 6.
It has been shown by Feige
[\ref{Fl}, \ref{Fu}] that
\[
(1+o(1))n\log n\le \E_v C\le
(1+o(1))\frac{4}{27}n^3,
\] 
and these bounds are tight.

In this paper, we show that for 
bounded-degree planar graphs, one has better bounds, namely,

\begin{theorem}\label{mainT}
Let $G=(V,E)$ be a finite connected planar graph with $n$ vertices
and maximal degree $M$.
Then for every vertex $v\in V$,
\[
c n (\log n)^2 < \E_v C< 6n^2\,,
\]
where $c$ is a positive constant depending only on $M$.
\end{theorem}

This generalizes a result  of Zuckerman [\ref{Z1}]
showing that $\min_{v}\E_{v}C \ge c n(\log n)^2$ for
bounded-degree trees on $n$ vertices.
If $G=\Z^d\cap[-m,m]^d$,
a finite portion of the $d$-dimensional integer lattice,
then $\E_v C$ is $\Theta(n^2)$ for $d=1$,
$\Theta(n(\log n)^2)$ for $d=2$ and $\Theta(n \log n)$ for $d \geq 3$
[\ref{A},\ref{Z2}].  Here, $n=(2m+1)^d=|V|$.
The cases $d=1$ and $d=2$ show that Theorem~\ref{mainT}
is tight (up to the constants).
The case $d=3$ shows that
the planarity assumption is necessary.

The upper bound in Theorem~\ref{mainT} is quite easy.
The lower bound will be based on Koebe's~[\ref{K}]
Circle Packing Theorem (CPT):

\begin{theorem}\label{CPT}
Let $G=(V,E)$ be a finite planar graph.
Then there is a disk packing $\bigl\{C_v:v\in V\bigr\}$
in $\RR^2$,
indexed by the vertices of $G$, such that
$C_v\cap C_u\neq\emptyset$ iff $\{v,u\}\in E$.
\end{theorem}

Koebe's proof relies on complex analysis, but recently
several new proofs have been discovered.
See, for example, [\ref{BrSc}] for a geometric, combinatorial
proof.

Some fascinating relations between the CPT
and analytic function theory have been studied
in the last decade.  Additionally, the CPT
became a tool for studying planar graphs
in general, and random walks on planar graphs
in particular [\ref{MT}, \ref{MP}, \ref{HS}, \ref{BS}, \ref{ST}].
In these applications, as well as here, the CPT is
useful because it endows the graph with a geometry
that is better, for many purposes, than the usual graph-metric.

\medskip
We conjecture that Theorem~\ref{mainT} holds with $c=c'/\log (M+2)$,
where $c'>0$ is a positive constant.
For example, this is true for trees, 
since in a tree one can easily find a set of at least
$n^{1/2}$ vertices with pairwise distances at least
$ \log n /\bigl(2 \log(M+2)\bigr)$.
As we shall see, this implies that the expected cover time
is bounded below by a constant times $n(\log n)^2/\log(M+2)$.

\section{Preliminaries}
For a simple random walk on the graph $G=(V,E)$ we define 
for every ordered pair $(u,v)$ of vertices, the {\em hitting time}
as $H(u,v) := \E_{u}T_{v}$.
The {\em commute time} is given by $C(u,v) := H(u,v)+H(v,u)$
and the {\em difference time} is given by 
\[
D(u,v) := H(u,v)-H(v,u).
\]
{}From the so called cyclic tour property of reversible Markov chains it
follows
that difference times are additive (see [\ref{CTW}]):
\begin{equation}\label{eA}
D(u,v)+D(v,w) = D(u,w).
\end{equation}

Commute times are closely related to {\em effective resistances\/} in 
electrical networks:
Regard each edge of $G$ as a unit resistor and define for each pair $(u,v)$
of vertices the effective resistance $R(u,v)$ between them as
$i^{-1}$ where $i$ is the current flowing into $v$ when
grounding $v$ and applying a 1 volt potential to $u$.
In mathematical terms, $R(u,v)$ can be defined as
\[
R(u,v):= \sup \frac{\bigl(f(v)-f(u)\bigr)^2}{{\cal D}(f)},
\]
where ${\cal D}(f)$ is the Dirichlet energy of $f$,
\[
{\cal D}(f):=
\sum_{\{a, b\}\in E}
\bigl(f(a)-f(b)\bigr)^2,
\]
and the $\sup$ is with respect to all $f:V\to\RR$ such that
${\cal D}(f)>0$.  (If $u$ and $v$ are in distinct components
of $G$, then $R(u,v)=\infty$.)
It is an immediate consequence from this definition
that when $G$ is a subgraph of another graph $G'$,
and $u,v$ are vertices in $G$, then the 
effective resistance between $u$ and $v$ in
$G'$ is bounded from above by the effective resistance
between them in $G$.
It is well known that resistances satisfy the triangle inequality
\begin{equation}\label{eTriang}
R(u,w)\le R(u,v)+R(v,w),
\end{equation}
which follows from the following useful formula from [\ref{CRRST}]:
\begin{equation} \label{eB}
C(u,v) = 2\left|E\right|R(u,v).
\end{equation}

There is also a formula from [\ref{Tet}] for $H(u,v)$ in terms of
resistances,
but it is more complicated:
\begin{equation}\label{eTet}
H(u,v) = \frac 12 \sum_{w\in V} d_w 
\bigl( R(u,v)+R(v,w)-R(u,w)\bigr),
\end{equation}
where $d_w$ is the degree of $w$.

The main lemma in the proof of Theorem~\ref{mainT} involves estimating
the resistances.
A combination of the above identities 
will then yield lower bounds for the hitting times.
We then need some way to estimate the cover time from the
hitting times.
For this {\em Matthews' method} [\ref{M}] will prove useful.

\begin{lemma} \label{lD}
Let $G(V,E)$ be a finite graph.  Then
\[\max_{v\in V}\E_v C \leq h_{n-1}\max\bigl\{H(u,v):v,u\in V\bigr\},\]
where $h_{k}$ denotes the harmonic series $\sum_{i=1}^{k}i^{-1}$.
Furthermore,
\[\min_{v\in V}\E_v C \geq h_{|V_0|-1}
\min\bigl\{H(u,v) : {u,v \in V_0,\, u \neq v}\bigr\}\]
holds for every subset $V_0\subset V$.
\end{lemma}

A proof can be found in [\ref{M}] or [\ref{AF}].
The proof relies on one ingenious trick,
namely, to assign a uniformly chosen random order to $V$ 
independent of the random walk.

\section{Proof of Theorem~\ref{mainT}}

We start with the easy proof of the upper bound.
It is a well known consequence of Euler's formula $|V|-|E|+|F|=2$
(see [\ref{D}, Theorem 4.2.7]) that the average degree
$\bar d$ in a finite planar graph is less than $6$.
By [\ref{AF}, Chapter 6, Theorem 1],
$\max_{v}\E_v C \leq \bar{d}n(n-1) < 6n^2$, which
gives the upper bound.

Let us now turn to the lower bound.
The main tool in the proof of Theorem \ref{mainT} is the following lemma:

\begin{lemma} \label{lR}
There exist positive constants $c=c(M)$ and $r=r(M)$
such that for every planar connected graph $G=(V,E)$ with 
maximum degree $M $ and every set of vertices
$W\subset V$ there is a subset $V'\subset W$
with $|V'| \geq |W|^{c}$ and 
$R(u,v) \geq r\log |W|$ for every $u\ne v$, $u,v\in V'$.
\end{lemma}

{\noindent\em Proof of Theorem \ref{mainT} from Lemma \ref{lR}.} 
The strategy is to convert the information Lemma \ref{lR}
gives about resistances to information about hitting times
$H(v,u)$, and then use the second part of Lemma~\ref{lD}.

Let $a\in V$ be some vertex, and let $\{v_1,v_2,\dots,v_n\}$
be an ordering of $V$ such that
$i\le j$ implies $D(a,v_i)\le D(a,v_j)$.
Then we have
\begin{equation}
i\le j\Rightarrow D(v_i,v_j)\ge 0
\label{eO}
\end{equation}
for all $i,j\in\{1,\dots,n\}$, by (\ref{eA}).
Let $k=[n/2]$, the largest integer in $[0,n/2]$.
We now consider several distinct cases.

{\medskip\noindent\bf Case 1:}
there are some $i<j$ in $\{1,2,\dots,n\}$ such that
$H(v_j,v_i)\ge n (\log n)^2/2$.
Observe that for all $v\in V$, we have
$\E_v C\ge\min\{H(v_j,v_i),H(v_i,v_j)\}$,
for the random walk starting at $v$ must either visit
$v_j$ before $v_i$ or visit $v_i$ before $v_j$.
Consequently, (\ref{eO}) completes the proof in this case.

{\medskip\noindent\bf Case 2:} $D(v_1,v_k)\ge n(\log n)^3$,
and Case 1 does not hold.
By (\ref{eO}) and (\ref{eA}), we then have 
$D(v_1,v_j)\ge n(\log n)^3$ for all $j\ge k$.
By (\ref{eTet}) and (\ref{eTriang}) we have
\begin{align*}
H(v_n,v_1)
&
\geq  \frac{1}{2} \sum_{j=k}^{n}
\Bigl(R(v_1,v_n)+R(v_1,v_j)-R(v_n,v_j)\Bigr)  
\\ &
= \frac{1}{4 |E|} \sum_{j=k}^{n}
\Bigl(C(v_1,v_n)+C(v_1,v_j)-C(v_n,v_j)\Bigr)  
\qquad\hbox{(by (\ref{eB}))} 
\\ &
> \frac{1}{12n} \sum_{j=k}^{n} \Bigl(C(v_1,v_n)+C(v_1,v_j)-D(v_j,v_n)-2
H(v_n,v_j)\Bigr),
\end{align*}
since $|E|< 3|V|$ for planar graphs
and $C(v_n,v_j)= 2H(v_n,v_j) + D(v_j,v_n)$.
Consequently, since Case 1 does not hold,
\begin{align*}
H(v_n,v_1)
&
> \frac{1}{12n}
\sum_{j=k}^{n} \Bigl(C(v_1,v_n)+C(v_1,v_j)-D(v_j,v_n)-n (\log n)^2\Bigr)  
\\ &
\geq \frac{1}{12n}
\sum_{j=k}^{n} \Bigl(D(v_1,v_n)+D(v_1,v_j)-D(v_j,v_n)-n (\log n)^2\Bigr)  
\\ &
= \frac{1}{12n} \sum_{j=k}^{n} \Bigl(2D(v_1,v_j)-n(\log n)^2\Bigr)
\geq \frac{1}{12} n(\log n)^3,
\end{align*}
for all sufficiently large $n$,
since we have $D(v_1,v_j)\ge n(\log n)^3$ for all $j\ge k$.
However, $H(v_n,v_1)\ge \frac{1}{12} n(\log n)^3$ brings
us back to Case 1.

{\medskip\noindent\bf Case 3:} $D(v_1,v_k)\le n(\log n)^3$.
Set $W=\{v_1,\dots,v_k\}$, and let $V'\subset W$ be as in Lemma \ref{lR}.
Let $m:=|V'|\ge n^c$, and let $i_1<i_2<\cdots <i_m$ be those
indices $i\in\{1,\dots,k\}$ such that $v_i\in V'$.
Set $s:=[\sqrt m]-1$.  Since
\[
\sum_{j=1}^{s-1} D(v_{i_{js}},v_{i_{(j+1)s}})
= D(v_{i_s},v_{i_{s^2}}) \le D(v_1,v_k) \le n(\log n)^3,
\]
there is some $t\in\{1,\dots,s-1\}$
such that $D(v_{i_{ts}},v_{i_{(t+1)s}})\le n(\log n)^3/(s-1) = o(n)$.
Set $V_0:=\{v_{i_{ts}},v_{i_{ts+1}},\dots,v_{i_{(t+1)s}}\}$.
Then $|V_0|\ge n^{c'}$ for some constant $c'>0$
and $D(u,w)\le o(n)$ for $u,w\in V_0$, if $n$ is large.
However, we have $C(u,w)=2 |E| R(u,w) \ge r n\log n$ 
for $u,w\in V_0$, since $V_0\subset V'$.
Because $2H(u,w)=C(u,w)- D(w,u)$, this gives
$H(u,w)\ge (r/3)n\log n$ for $u,w\in V_0$, provided that $n$ is
large.
Now the second part of Lemma \ref{lD} completes the proof.  \hfill{$\Box$}

\medskip\noindent 
{\bf Remark.} The recent preprint by Kahn et.\ al.\ [\ref{KKLV}]
gives an estimate (Prop.~1.2 and Thm.~1.3) of the
expected cover time in terms of the commute times. 
This result could be used to simplify the above argument
(but was not available at the time of writing of the
first draft of the current paper).

\bigskip

{\noindent\em Proof of Lemma \ref{lR}.} 
We first consider the case where $G$ is a triangulation
of the sphere.  This means that $G$ is a graph embedded in $S^2$ 
with the property that every connected component
of $S^2\setminus G$ has precisely $3$ edges of $G$ as its
boundary. 

The Circle Packing Theorem 
implies the existance of a disk packing $\bigl\{C_v:v\in V\bigr\}$
indexed by the vertices of $G$, such that
each $C_v$ is a closed round disk in $\RR^2$ and
$C_v\cap C_u\neq\emptyset$ iff $\{v,u\}\in E$.
Moreover (by normalizing by a M\"obius transformation),
we assume with no loss of generality that the
outer three disks in the packing all have radius $1$.

The Ring Lemma from [\ref{RS}] 
implies that there is a constant
$A = A(M)$ such that
\begin{equation}\label{ratio}
\{v,u\}\in E\Rightarrow r_v < A r_u,
\end{equation}
where $r_v$ denotes the radius of $C_v$.
It then follows that there is another constant $A' := A'(M)>0$
such that
\begin{equation}\label{dist}
\{v,u\}\notin E\Rightarrow
\dist(C_v,C_u)\ge A' r_u,
\end{equation}
where $\dist(C_v,C_u):=\inf\{|p-q|:p\in C_v,\,q\in C_u\}$,
because the disks around $C_u$ separate $C_v$ from $C_u$,
since $G$ is assumed to be a triangulation.

Most important for us is the following lower bound for the resistance
\begin{equation}
R(w,u)\ge A''\log\bigl( \dist(C_w,C_u)/r_u\bigr),
\label{er}
\end{equation}
for some constant $A'' := A''(M) >0$.
Similar estimates appear in [\ref{HS}] and in [\ref{BS}].
For completeness, we include a quick proof here. For each $v\in V$ let
$z_v$ be the center of the disk $C_v$.   
Set $a:=\log r_u$ and $b:= \log |z_w-z_u|$.
Consider $F(z):= \log(z-z_u)$
as a map from $\C\setminus\{z_u\}=\RR^2\setminus\{z_u\}$ to the cylinder
$\RR+i (\RR/2\pi\Z)$. 
Set $f(v):= \min\{\Re F(z_v),b\}$, for $v\neq u$
and $f(u):= \log r_u=a$.
The inequality $\dist(C_u,C_v)\ge A' r_v$
implies that $\area\bigl(F(C_v)\bigr)/ \diam\bigl(F(C_v)\bigr)^2$
is bounded above and below by positive constants.
For neighbors $v_1$ and $v_2$ we have
\[
|f(v_1)-f(v_2)|\le \diam\bigl(F(C_{v_1})\bigr)+
\diam\bigl(F(C_{v_2})\bigr)
\le O(1) \diam\bigl(F(C_{v_1})\bigr)
\]
and $f(v_1)-f(v_2)=0$ unless
$\dist(C_u,C_{v_1}\cup C_{v_2})\le |z_w-z_u|-r_u$.
Consequently, ${\cal D}(f)\le O(1)\sum_v \diam\bigl(F(C_v)\bigr)^2$,
where the sum extends over all $v\neq u$ such that $F(C_v)$
intersects the cylinder $[a,b] + i (\RR/2\pi\Z)$.
All these sets $F(C_v)$ are contained in the cylinder
$\bigl[a,b+O(1)\bigr]+i(\RR/2 \pi\Z)$, 
and their interiors are disjoint.
Since the area of each $F(C_v)$ is proportional
to the square of its diameter, we find that
\[{\cal D}(f)\le O(1)\area\bigl([a+O(1)]+i(\RR/2 \pi\Z)\bigr)
= O(1)(b-a+1).
\]
The inequality (\ref{er}) now follows from the definition
of the effective resistance.

Fix a small $s>0$ (which will be specified later),
and set $n=|W|$.
For $j\in\Z$, let
\[
W_j:=\bigl\{v\in W: r_v\in(n^{s(j-1)},n^{sj}]\bigr\}.
\]
Then $W=\bigcup_{j\in\Z} W_j$.
For $n$ so large that $n^s \geq A$ we have by (\ref{ratio}) that
if $u \in W_{j}$, $v \in W_{k}$ and $k-j \geq 2$ then
$\{u,v\} \not \in E$, and by (\ref{dist}) and (\ref{er}),
$R(u,v) \geq A''\log(\dist(C_{u},C_{v})/r_u) 
\geq A'' \log (A' r_{v}/r_u)
\geq \frac{1}{2}A'' s\log n$,
when $n$ is large.

Now either
$\left|\bigcup_{j\ \hbox{\rm\scriptsize odd}}W_{j}\right| \geq n/2$
or 
$\left|\bigcup_{j\ \hbox{\rm\scriptsize even}}W_{j}\right| \geq n/2$.
Let us assume the latter case, noting that the former is treated
similarly.

For each even $j$, let $Z_{j}$ be a maximal subset of 
vertices of $W_{j}$ such that 
\[
u,v \in Z_{j},\ u\ne v \quad\Rightarrow \quad
\dist(C_{u},C_{v}) \geq n^{s(j+1)},
\]
and note that by the definition of $W_j$ and (\ref{er}),
$R(u,v) \geq A''\log(n^{s(j+1)}/n^{sj}) = A'' s\log n$
for all $u,v \in Z_{j}$, $u\neq v$.
Since for any $v \in W_{j}$
the disk of radius $3n^{s(j+1)}$ centered at $z_v$, the center of 
$C_{v}$, does not contain more than
$\bigl(3n^{s(j+1)}/n^{s(j-1)}\bigr)^2 = 9 n^{4s}$  
disks
$C_{u}$ with $u \in W_{j}$, it follows that $|Z_{j}| \geq n^{-4s}|W_{j}|/9$.
Now put $V' = \bigcup_{j \ \hbox{\rm\scriptsize even}}Z_{j}$.
Then $|V'| \geq n^{1-5s}$ for $n$ large enough and when
$v\ne v'$ are in $V'$
we have $R(u,v) \geq \frac{1}{2}A''s\log n$.
The result for $G$ a triangulation
of $S^2$ follows by choosing $s=1/6$, say.

Now consider the case where $G$ is not a triangulation
of $S^2$.
It is easy then to construct a triangulation $T$ of the sphere
with maximum degree at most $3M$ which contains $G$
as a subgraph.  
The effective resistence $R_G(u,v)$ in $G$ between two vertices 
$u,v$ in $G$ is at least $R_T(u,v)$, their effective resistance
in $T$.  Consequently, this case follows from the previous.

\hfill{$\Box$}

\begin{center}
{\sc References}
\end{center}

\begin{enumerate}
\item {\sc D. Aldous}, On the Time Taken by Random Walks on Finite Groups to
Visit Every State, {\em Wahrsch.\ Verw.\ Gebeite} {\bf 62} (1983), 361--374.
\label{A}
\item {\sc D. Aldous  and J. A. Fill}, ``Reversible Markov Chains and Random
Walks on Graphs'' (book draft), October 1999. \hfill\break
{\texttt http://www.stat.berkeley.edu/\string~aldous/book.html} \label{AF}
\item {\sc I. Benjamini and O. Schramm}, Harmonic Functions on Planar and
Almost Planar Graphs and Manifolds, Via Circle Packings, {\em Invent.\
Math.} {\bf 126} (1996), 565--587. \label{BS}
\item {\sc G. R. Brightwell and E. R. Scheinerman}
Representations of planar graphs,
{\em SIAM J.\ Discrete Math.} {\bf 6} (1993), 214--229.
\label{BrSc}
\item {\sc A. Chandra, P. Raghavan, W. Russo, R. Smolensky and P. Tiwari},
The electrical
resistance of a graph captures its commute and cover times, {\em Comput.\
Complexity} {\bf 6} (1996/97), 312--340. \label{CRRST}
\item {\sc D. Coppersmith, P. Tetali and P. Winkler},
Collisions Among Random Walks on a Graph, {\em SIAM J.\ Discrete
Math.} {\bf 6} (1993), 363--374. \label{CTW}
\item {\sc R. Diestel}, ``Graph Theory,'' Springer, New York, 1997.
\label{D}
\item {\sc U. Feige}, A Tight Upper Bound on the Cover Time for Random
Walks on Graphs, {\em Random Struct.\ Alg.} {\bf 6} (1995), 51--54.
\label{Fu}
\item {\sc U. Feige}, A Tight Lower Bound on the Cover Time for Random
Walks on Graphs, {\em Random Struct.\ Alg.} {\bf 6} (1995), 433--438.
\label{Fl}
\item {\sc Z. He and O. Schramm}, Hyperbolic and Parabolic Packings, {\em
Discrete Comput.\ Geom.} {\bf 14} (1995), 123--149. \label{HS}
\item {\sc J. Kahn, J. H. Kim, L. Lov\'{a}sz and V. H. Vu}, The cover time, the blanket time, and the Matthews bound, Preprint (2000). \label{KKLV} 
\item {\sc P. Koebe}, Kontaktprobleme der konformen abbildung, {\em Berichte
Verhande.\ S\"{a}chs.\ Akad.\ Wiss.\ Leipzig, Math.-Phys.\ Klasse} {\bf 88}
(1936) 141--164. \label{K}
\item {\sc S. Malitz and A. Papakostas},
{On the angular resolution of planar graphs},  
{\em SIAM J.\ Discrete Math.},  
{\bf 7}  
(1994),
{172--183}. \label{MP}
\item {\sc P. Matthews}, Covering Problems for Brownian Motion on 
Spheres, {\em Ann.\ Probab.} {\bf 16} (1988), 189--199. \label{M}
\item {\sc G. L. Miller and W. Thurston}, Separators in Two and
three dimensions, in Proc.\ of the 22th Annual ACM Symposium
on Theory of Computing, 300--309, Baltimore, (1990), ACM. \label{MT}
\item {\sc B. Rodin and D. Sullivan},
The convergence of circle packings to the Riemann mapping,
 {\em J.\ Differential Geom.} {\bf 26} (1987),  349--360. \label{RS}
\item {\sc D.A.~Spielman and S.-H.~Teng},
Spectral partitioning works: planar
graphs and finite element meshes, {\em Proc.\ 37th Ann.\ Symp.\
Found.\ of Comp.\ Sci.}, IEEE (1996), 96--105. \label{ST}
\item {\sc P. Tetali},
Random walks and the effective resistance of networks, 
{\em J.\ Theoret.\ Probab.} {\bf 4} (1991) 101--109. \label{Tet}
\item {\sc D. Zuckerman}, Covering Times of Random Walks on
Bounded-Degree Trees and Other Graphs, {\em J.\ Theor.\ Probab.} {\bf 2}
(1989), 147--158. \label{Z1}
\item {\sc D. Zuckerman}, A technique for lower bounding the cover 
time, {\em SIAM J.\ Disc.\ Math.} {\bf 5} (1992), 81--87. \label{Z2}
\end{enumerate}

\bigskip
\hbox{
\vtop{\hsize=2.5in
\small\noindent
Johan Jonasson \\
Dept.\ of Mathematics \\
Chalmers University of Technology \\
S-412 96 Sweden \\
Phone: +46 31 772 35 46 \\
\texttt{jonasson@math.chalmers.se}
}
\hskip 1in
\vtop{\hsize=2.5in
\small\noindent
Oded Schramm\\
Microsoft Corporation,\\
One Microsoft Way,\\
Redmond, WA 98052; USA\\
\texttt{schramm@Microsoft.com}
}
}
\end{document}